
\input gtmacros
\input amsnames
\input amstex

%
\catcode`\@=12        
\input gtmonout
\volumenumber{2}
\volumeyear{1999}
\volumename{Proceedings of the Kirbyfest}
\pagenumbers{87}{102}
\papernumber{4}
\received{5 March 1999}\revised{24 June 1999}
\published{21 October 1999}

%
\let\\\par
\def\topmatter{\relax}

\let\gttitle\title
\def\title#1\endtitle{\gttitle{#1}}
\let\gtauthor\author
\def\author#1\endauthor{\gtauthor{#1}}
\let\gtaddress\address
\def\address#1\endaddress{\gtaddress{#1}}
\def\affil#1\endaffil{\gtaddress{#1}}
\let\gtemail\email
\def\email#1\endemail{\gtemail{#1}}
\def\subjclass#1\endsubjclass{\primaryclass{#1}}
\let\gtkeywords\keywords
\def\keywords#1\endkeywords{\gtkeywords{#1}}
\def\heading#1\endheading{{\def\S##1{\relax}\def\\{\relax\ignorespaces}
    \section{#1}}}
\def\head#1\endhead{\heading#1\endheading}

\def\subhead#1\endsubhead{\sh{#1}}
\def\subsubhead#1\endsubsubhead{\sh{#1}}
\def\specialhead#1\endspecialhead{\sh{#1}}

\def\qed{\ifmmode\quad\sq\else\hbox{}\hfill$\sq$\par\goodbreak\rm\fi}  
\def\proclaim#1{\rk{#1}\sl\ignorespaces}
\def\endproclaim{\rm\ppar}
\def\cite#1{[#1]}
\newcount\itemnumber

\let\itemold\item
\def\item{\itemold{{\rm(\number\itemnumber)}}%
\global\advance\itemnumber by 1\ignorespaces}
\def\S{section~\ignorespaces}  
\def\date#1\enddate{\relax}
\def\thanks#1\endthanks{\relax}   
\def\dedicatory#1\enddedicatory{\relax}  
\let\footnote\plainfootnote

\def\Refs{\ppar{\large\bf References}\ppar\bgroup\leftskip=25pt
\frenchspacing\parskip=3pt plus2pt\small}       
\def\endRefs{\egroup}
\def\widestnumber#1#2{\relax}
\def\endrefitem{}
\def\refdef#1#2#3{\def#1{\leavevmode\unskip\endrefitem#2\def\endrefitem{#3}}}
\def\ref{\par}
\def\endref{\endrefitem\par\def\endrefitem{}}
\refdef\key{\noindent\llap\bgroup[}{]\ \ \egroup}
\refdef\no{\noindent\llap\bgroup[}{]\ \ \egroup}
\refdef\by{\bf}{\rm, }
\refdef\manyby{\bf}{\rm, }
\refdef\paper{\it}{\rm, }
\refdef\book{\it}{\rm, }
\refdef\jour{}{ }
\refdef\vol{}{ }
\refdef\yr{$(}{)$ }
\refdef\ed{(}{ Editor) }
\refdef\publ{}{ }
\refdef\inbook{from: ``}{'', }
\refdef\pages{}{ }
\refdef\page{}{ }
\refdef\paperinfo{}{ }
\refdef\bookinfo{}{ }
\refdef\publaddr{}{ }
\refdef\moreref{}{ }
\refdef\eds{(}{ Editors)}
\refdef\bysame{\hbox to 3 em{\hrulefill}\thinspace,}{ }
\refdef\toappear{(to appear)}{ }
\refdef\issue{no.\ }{ }
\TagsOnRight

\topmatter
\title Topological Field Theories and formulae\\of Casson and
Meng--Taubes
\endtitle
\shorttitle{Topological Field Theories}
\asciititle{Topological Field Theories and formulae of Casson and
Meng-Taubes}

\author S\thinspace K Donaldson
\endauthor
\asciiauthors{S K Donaldson}

\affil Department of Mathematics\\Imperial College, 
London SW7 2BZ, UK\endaffil
\email s.donaldson@ic.ac.uk \endemail

\abstract The goal of this paper is to give a new proof of a theorem of
Meng and Taubes [9] that identifies the Seiberg--Witten invariants of
3--manifolds with Milnor torsion. The point of view here will be
that of topological quantum field theory.  In particular, we relate the
Seiberg-Witten equations on a 3--manifold with the Abelian vortex equations
on a Riemann surface. These techniques also give a new
proof of the surgery formula for the Casson invariant,
interpreted as an invariant of a  homology $S^2\times S^1$.
\endabstract
\asciiabstract{The goal of this paper is to give a new proof of a theorem of
Meng and Taubes that identifies the Seiberg-Witten invariants of
3-manifolds with Milnor torsion. The point of view here will be
that of topological quantum field theory.  In particular, we relate the
Seiberg-Witten equations on a 3-manifold with the Abelian vortex equations
on a Riemann surface. These techniques also give a new
proof of the surgery formula for the Casson invariant,
interpreted as an invariant of a  homology S^2 x S^1.}

\primaryclass{57R57}
\secondaryclass{57M25, 57N10, 58D29}

\keywords
Seiberg--Witten invariant, Casson invariant, Alexander polynomial,
Milnor torision, topological quantum field theory, moduli space, vortex
equation
\endkeywords
\asciikeywords{Seiberg-Witten invariant, Casson invariant, Alexander 
polynomial, Milnor torsion, topological quantum field theory, moduli
space, vortex equation}

\maketitle
\document
\NoBlackBoxes

\section{Introduction}

 In 1985, Casson introduced his renowned invariant of homology
3--spheres, together with  a surgery formula for the change in the
invariant when the manifold is changed by a surgery [1]. The latter
can also be seen as a formula for a Casson-type invariant of a
homology $S^{1}\times S^{2}$. For a manifold $Y^{3}$ of this kind we
\lq\lq count'', in the sense of Fredholm differential topology---as
in the work of Taubes  [11]---the flat connections on a non-trivial
$SO(3)$ bundle over $Y$ to obtain a number
$C(Y)$. Casson's formula is then
$$   C(Y) = \sum_{i>0} i^{2} a_{i}, \tag1$$ where the integers
$a_{i}$ are the coefficients of the normalised Alexander polynomial
$\Delta_{Y}(t) = a_{0} + \sum_{i>0} a_{i} (t^{i} + t^{-i})$. (The
interpretation of Casson's surgery formula in this way is implicit
in the Floer's exact sequence [4] relating the Floer homologies of
manifolds differing by surgeries.)

In 1995, Meng and Taubes [9] considered the dimensionally-reduced
Seiberg--Witten equations over homology $S^{1}\times S^{2}$'s. Here
one has an integer parameter $d>0$ and invariants $SW_{d}(Y)$ which
count the solutions of the equations on a complex line bundle of
Chern class $d$ over $Y$. The Meng--Taubes formula can be written
$$     SW_{d}(Y) = a_{d+1}+ 2 a_{d+2} + 3 a_{d+3} +\dots    .\tag2$$

The purpose of this article is to describe how the ideas of
Topological  Field Theory can be applied to derive these formulae in
ways that are, at least superficially, different from the original
proofs. The author has developed this point of view intermittently
over the last 10 years. The  Casson case was worked out  in
discussions with M Furuta around 1989, partly inspired by lectures
of Segal on Quantum Field Theory (and some of the constructions
described in Section 5 are due to Furuta).  The Seiberg--Witten case
was worked out in 1995, partly motivated by the
 programme of Pidstragatch and Tyurin and work of Thaddeus, which we
mention again at the end of the article. Over the years a number of
authors have developed ideas which come very close to the
 point of view we take here: we mention particularly the work of
Frohman and Nicas [5]
 on the one hand and Hutchings and Lee [7] on the other. But we hope
that the approach here may add  something which does not quite
appear in the literature and which may be worth recording; without
staking any particular claims to originality. To keep the exposition
short and simple we take the liberty of treating the
differential-geometric  and analytical foundations of our
arguments---which we believe can all be
 pieced together from the literature---very sketchily, and
concentrate on the formal aspects. (As a  particular instance of
this sketchiness, we will not go into the question of {\it
orientations}, so our formulae will really be derived
 up to overall sign $\pm1$.)

\section{Topological field theories}

We will consider theories which satisfy the Topological Field Theory
axioms [2] in a restricted range of cases. In essence we want to
associate to each  closed, connected, oriented surface $\Sigma$ a
vector space $V(\Sigma)$ and to a cobordism
$W$ between two such---$\Sigma_{0},\Sigma_{1}$---a linear map
$$    \rho_{W}\co V(\Sigma_{0})\rightarrow V(\Sigma_{1}), $$ satisfying
the familiar  composition rule when one glues together two
cobordisms. In our examples the spaces will come with natural
gradings,
 although  we will need only keep track of the $\bold{Z}/2$--grading.
In the case when $\Sigma_{0}, \Sigma_{1}$ are diffeomorphic, and we
choose a definite diffeomorphism  between them, we may form a closed
$3$--manifold
$\overline{W}$. We consider theories where we have numerical
invariants for such closed manifolds $n(\overline{W})$ and the other
gluing axiom we require is that
$$    n(\overline{W}) = \text{Tr}_{s}( \rho_{W}), $$ where
$\rho_{W}$ is regarded as an endomorphism of $V(\Sigma_{0})$, using
the fixed  identification of  $\Sigma_{0},\Sigma_{1}$ and
$\text{Tr}_{s}$ denotes a
\lq\lq supertrace'', with signs induced by the grading.

We will now recall how  formal structures of this kind
 arise in gauge theory. We begin with the set-up which will
correspond to the Casson formula. For each  surface $\Sigma$ we
construct the moduli space
 $M(\Sigma)$ of flat connections on a non-trivial $SO(3)$ bundle
over $\Sigma$. This is an orbifold of dimension $6g(\Sigma)-6$. Our
vector space $V^{C}(\Sigma)$ will then be the homology of this
moduli space.  A technical variant of this construction is to
consider the moduli space of projectively flat connections on a
$U(2)$ bundle of odd-degree over $\Sigma$, with fixed central
curvature. The moduli space of these connections is a smooth
manifold $\tilde{M}(\Sigma)$ covering (in the orbifold sense)
 $M(\Sigma)$. More precisely, $M(\Sigma)$ is the quotient of
$\tilde{M}(\Sigma)$ by a natural action of $H^{1}(\Sigma;
\bold{Z}/2)$.  It is somewhat easier to work in this
$U(2)$ framework, since one does not run into the notorious
difficulties with reducible connections, but the language becomes
more tortuous.
 In the end it does not really matter how one proceeds since we will
only be concerned with the real homology and it is known that
 the real homology groups of $M(\Sigma), \tilde{M}(\Sigma)$ are
isomorphic. Equivalently the action of $H^{1}(\Sigma; \bold{Z}/2)$
on the real homology is trivial  [3]. Likewise, we know that
Poincar\'e duality and intersection theory goes over for the real
homology of an orbifold. To simplify our exposition we will work
with the
$SO(3)$ moduli spaces and  ignore the fact that they are not quite
manifolds.

Turning to $3$--manifolds, for a closed, oriented,  $3$--manifold $Y$
with the
 homology  of $S^{1}\times S^{2}$  we consider of course the
invariant
 $C(Y)$ obtained by counting the flat
$SO(3)$ connections on a non-trivial bundle. Note that the
topological hypotheses  mean that there is a unique non-trivial
$SO(3)$ bundle over $Y$, and there are no  reducible flat
connections on this bundle. We could extend the theory to certain
 other $3$--manifolds, but we would need to specify the relevant
bundle and avoid reducibles, so for simplicity we will stick to the
class of homology $S^{1}\times S^{2}$'s. Now let $W$ be a cobordism
from
$\Sigma_{0}$ to $\Sigma_{1}$. We will restrict attention to the
cobordisms
$W$ such that $H_{1}(\Sigma_{0})\oplus H_{1}(\Sigma_{1}) $ maps {\it
onto}
$H_{1}(W)$ under the homorphisms induced by inclusion, or
equivalently that
$$    H_{1}(W, \partial W)= \bold{Z}. \tag3$$ All the manifolds we
encounter---occurring the decomposition of a homology $S^{1}\times
S^{2}$---will satisfy this condition. The main idea is this: we
consider a moduli space
$M(W)$ of flat connections over $W$, with a restriction map
$$ r\co M(W) \rightarrow M(\Sigma_{0})\times M(\Sigma_{1}). $$
 If $M(W)$ is a compact, oriented manifold it carries a fundamental
homology class which can be pushed forward to
 $r_{*}[M(W)]\in H_{*}(M(\Sigma_{0}\times M(\Sigma_{1})$. Then
Poincare duality and the Kunneth theorem give an identification
$$\eqalign{H_{*}(M(\Sigma_{0})\times M(\Sigma_{1}))& \cong
 \text{Hom}\ ( H_{*}(M(\Sigma_{0}), H_{*}(M(\Sigma_{1}))\cr&=
\text{Hom}\ (V^{C}(\Sigma_{0}), V^{C}(\Sigma_{1})),} $$ and we define
$\rho_{W}$ to be map obtained in this way from $r_{*}[M(W)]$.  Of
course all of this needs to be interpreted with the usual caveats of
the subject: we need to have a framework in which the equations
defining
  the moduli space   fit into a \lq\lq Fredholm package'', with a
well-defined index or formal dimension, and we may need to perturb
the equations to obtain transversality. There are two ways in which
this can be done. In one approach, as in [11], one sets up an
elliptic boundary-value problem on the manifold-with-boundary $W$.
In the other one adds semi-infinite cylinders to the ends to
construct a complete manifold $\hat{W}$. Then the set-up can be seen
as a dimension-reduced version of the $(3+1)$--dimensional Floer
theory, regarding flat connections over
$\hat{W}$ as instantons over $\hat{W}\times S^{1}$.  In either case
one finds that the relevant index or formal dimension is $3
\chi(\partial W)$. The functorial property for compositions of
cobordisms, and the relation with the numerical invariant
$C(\overline{W})$ follow from the usual \lq\lq gluing theory''
arguments. Notice that any homology--$S^{1}\times S^{2}$ can be
realised as a manifold
$\overline{W}$, where $W$ satisfies (3).
 \define \sym {\text{Sym}\ }

The Seiberg--Witten case follows a similar pattern. Here it seems
that one has to work in the tubular-end framework. We choose a
Riemannian metric and spin structure on the manifold $\hat{W}$ and a
line bundle $L\rightarrow \hat{W}$ with
$$    \langle c_{1}(L), [\Sigma_{0}]\rangle  = \langle c_{1}(L),
\Sigma_{1}\rangle = d. $$ Then we have Seiberg--Witten equations for
a pair consisting of a connection and an $L$--valued spinor field
over $\hat{W}$. If $d>0$ we do not run into  difficulties with
reducible solutions. To find the appropriate moduli spaces
associated with the  boundary components we look at the
translation-invariant solutions of these equations, which are the
potential limits of finite-energy solutions over the cylindrical
ends. These translation-invariant solutions over a tube
$\Sigma\times\bold{R}$ are the solutions of a {\it vortex equation}
over
$\Sigma$,  [8]. That is we fix a spin structure and metric on
$\Sigma$ and  look for pairs $(a,\psi)$ consisting of a connection
$a$ on
$L\rightarrow\Sigma$ and a section $\psi$ of $L^{*}\otimes
S^{+}(\Sigma)$ such that
  $$\align D_{a}\psi&= 0\\
           *F(a)  &= \vert \psi\vert^{2}. \endalign$$ Recall that
over the surface $\Sigma$ the spinor bundle $S^{+}$ can be identified
with the square root $K_{\Sigma}^{1/2}$ and the Dirac operator
$D_{a}$ can be identified with the Cauchy--Riemann operator, so that
$\psi$ becomes a holomorphic section of a holomorphic line bundle
$L^{*}\otimes K^{1/2}$ of degree $k= g-1-d$. Thus the zero set of
$\psi$ is a positive divisor of degree $k$, or equivalently an
element of the symmetric product $\sym^{k}(\Sigma)$. The basic
result we need  is that this sets up a 1--1 correspondence between
the moduli space of solutions of the vortex equation and the
symmetric product, provided $d>0$. (This is interpreted as saying
that the moduli space is empty if $d>g-1$.)
 We see then that the vector space $V_{d}(\Sigma)$
 we should associate to $\Sigma$ in this theory, for a given
parameter $d$, is the real homology
$H_{*}(\sym^{k}(\Sigma))$, where $k=g-1-d$. With this said, the
general construction goes through much as before, appealing to the
literature on the Seiberg--Witten--Floer theory for the relevant
gluing and transversality results.

\section{Homology of the moduli spaces}

     An important  part of the general topological field theory
package is the fact that the vector space $V(\Sigma)$ carries a
natural representation of the mapping class group of the surface
$\Sigma$. (More precisely, in the Seiberg--Witten case we should
consider finite-index subgroups of  diffeomorphisms which preserve a
spin structure.)
 The essence of our approach to the
$3$--dimensional invariants is to obtain a good understanding of
these  representations. We do the Seiberg--Witten case first.

\sh{(i) Homology of symmetric products}

                  For any space $X$ the real (or rational) homology
of the symmetric products $\sym^{k}(X)$ is easy to describe. We
begin with the homology of the $k$--fold product $X\times \dots
\times X$, which is just the tensor product of $k$ copies of
$H_{*}(X)$, and then take the  invariant part under the action of
the permutation group. In the case of a surface this gives
$$  H_{*}(\sym^{k}(\Sigma)) = (\Lambda^{k}\otimes s^{0}) \oplus
(\Lambda^{k-1}\otimes s^{1})\oplus \dots \oplus s^{k}, \tag4$$ where
$\Lambda^{j}=\Lambda^{j}(H_{1}(\Sigma)) $ and $s^{j}$  is the jth.
symmetric power of $H_{0}(\Sigma)\oplus H_{2}(\Sigma)$. So if we
write
$u$ for the fundamental class of $\Sigma$,
$$  s^{j} = \langle 1, u, \dots, u^{j}\rangle, $$ and
$\text{dim}(s^{j})= j+1$.

It will be convenient to reorganise the grading of our exterior
powers by writing
$$ \Lambda_{(i)} = \Lambda^{g-i},    $$ so we have isomorphisms
$\Lambda_{(i)}\cong \Lambda_{(-i)}$ induced by  the intersection
form on $H_{1}(\Sigma)$. Thus, remembering that
$k= g-1-d$, our formula (4) becomes
$$   V_{d}(\Sigma)= \Lambda_{(d+1)} \oplus \langle 1,u\rangle
\otimes \Lambda_{(d+2)} \oplus \langle 1,u,u^{2}\rangle \otimes
\Lambda_{(d+3)}\dots, $$ or, in shorthand,
$$   V_{d}(\Sigma) =
 \Lambda_{(d+1)} + 2 \Lambda_{(d+2)} + 3 \Lambda_{(d+3)}\dots
.\tag5$$ Here the  reader will immediately see the shape of the
 Meng--Taubes formula (2)  appearing. Tracing through the
constructions one verifies easily
 enough that these formulae describe the vector spaces
 $V_{d}(\Sigma)$ as representations of the mapping class group.
Notice that (although we will not need this) we can also bring in
the grading of the homology groups of the symmetric products. Again
it is best to reorganise the standard grading by writing:
$$   V_{d}^{(i)}(\Sigma)= H_{k-i}( \sym^{k}(\Sigma)),
     \tag6$$ so $V_{d}^{(0)}$ is the middle-dimensional homology. A
convenient notation is then to write
$$   \sum V_{d}^{(i)}(\Sigma) t^{i} = \Lambda_{(d+1)}+ (t+t^{-1})
\Lambda_{(d+2)} +  (t^{2} + 1 + t^{-2}) \Lambda_{(d+3)} +\dots, \tag
7$$ where we mean that the vector space appearing as the
coefficients of $t^{i}$ on the left hand side is isomorphic to the
corresponding  direct sum on the right hand side.

\sh{(ii) Homology of moduli spaces of flat connections}

 We begin with the formulae, obtained by Atiyah and Bott and other
authors using a variety of methods, for the Betti numbers of the
moduli space $M(\Sigma)$. The Poincar\'e polynomial is
$$\eqalign{P_{t}(M(\Sigma))& =\sum \text{dim}\ (H_{j}(M(\Sigma))) t^{j}\cr
&= \frac{1}{(1-t^{2})(1-t^{4})} \left( (1+t^{3})^{2g} -
t^{2g}(1+t)^{2g}\right).\cr} $$ Notice that the coefficients of $t^{j}$
in the expansion of $(1+t)^{2g}$ are the dimensions of the exterior
powers
$\Lambda^{j}$.  Thus we have
$$\sum_{j} \text{dim} (H_{j}(M(\Sigma))) t^{j} =
\frac{1}{(1-t^{2})(1-t^{4})}
 \sum \text{dim} (\Lambda^{i})(t^{3i}- t^{2g+i}).   $$ This formula
strongly suggests that there should be a natural (ie mapping class
group equivariant) identification of the homology groups themselves,
which we can write in the notation introduced above  as
$$  \sum_{i} H_{i}(M(\Sigma)) t^{i}  =
 \frac{1}{(1-t^{2})(1-t^{4})} \sum_{j} \Lambda^{j}
 (t^{3j}- t^{2g+j}), $$
 Analysis of the argument of Atiyah and Bott shows that this is
indeed the case: the starting point is the fact that the homology of
the moduli space injects into the homology of the space of all
irreducible connections where the action of the mapping class group
is clear. The formulae are again simpler if we use the reorganised
grading $\Lambda_{(j)}$. Similarly, we write
$V^{(C,i)}(\Sigma)$ for the homology group $H_{(3g-3-i)}(M(\Sigma))$.
 Then we have
$$\eqalignno{  \sum_{i}& V^{(C,i)}(\Sigma) t^{i}\cr &= t^{-(3g-3)}
\frac{1}{(1-t^{2})(1-t^{4})}
\sum_{j\geq 0} \Lambda_{(j)} (t^{3g+3j} + t^{3g-3j}-
t^{3g+j}-t^{3g-j}),\cr
    &=  \sum_{j>0} \frac{ (t^{2j}-t^{-2j})(t^{j}-t^{-j})}
{(t^{2}-t^{-2})(t-t^{-1})} \Lambda_{(j)}(\Sigma).&(8)\cr}
$$ Writing
$$   \frac{(t^{2j}-t^{-2j})(t^{j}-t^{-j})}{(t^{2}-t^{-2})(t-t^{-1})}
=  (t^{2j-2} + t^{2j-4}+\dots
 + t^{4-2j}+ t^{2-2j})(t^{j-1}+\dots + t^{1-j}), $$ we can evaluate
at $t=1$ to get a formula for the total homology
$$ V^{C}(\Sigma)=  H_{*}(M(\Sigma)) = \sum j^{2} \Lambda_{(j)}.
\tag9$$

\section{Scheme of proof}

We begin by bringing the Alexander polynomial into the picture, and
for this we recall some multilinear algebra, which one can find in
[10].  Suppose
$U_{0}, U_{1}$ are finite-dimensional vector spaces and $\Gamma$ is
a linear subspace of $U_{0}\oplus U_{1}$. Then, up to a scalar
factor, $\Gamma$ defines a Plucker point $\vert \Gamma\vert$ in the
exterior algebra $\Lambda^{*}(U_{0}\oplus U_{1})$. In turn, up to a
scalar ambiguity, elements of this exterior algebra can be viewed as
linear maps from $\Lambda^{*}(U_{0})$ to $\Lambda^{*}(U_{1})$, so we
have
$$ \vert \Gamma\vert \co  \Lambda^{*}(U_{0})\rightarrow
\Lambda^{*}(U_{1}), $$ defined up to a scalar. More precisely,
$\vert \Gamma\vert$ is defined once one chooses a trivialisation of
the line
$\Lambda^{\max}(\Gamma)\otimes \Lambda^{\max}(U_{0})^{*}$. In the
case when $\Gamma$ is the graph of a linear map $f\co U_{0}\rightarrow
U_{1}$ there is a natural trivialisation of this line and the map
$\vert \Gamma\vert$ is the usual map $\Lambda^{*}(f)$ induced on the
exterior powers. In general the construction satisfies a composition
rule, for the case when one has another subspace
$\Gamma'\subset U_{1}\oplus U_{2}$, provided the subspaces satisfy a
transversality condition. One wants the sum of the projection maps
from $\Gamma\oplus \Gamma'$ to $U_{1}$ to be surjective. Then if one
defines $\Gamma''\subset U_{0}\oplus U_{2}$
 to be set of pairs $(u_{0},u_{2})$ for which there exists a
$u_{1}\in U_{1}$ with
  $(u_{0},u_{1})\in \Gamma, (u_{1},u_{2})\in \Gamma'$, one has
$\vert \Gamma''\vert = \vert \Gamma'\vert\circ \vert \Gamma \vert$,
if one uses  an appropriate rule for normalising the scalar
ambiguities. We apply these ideas in the case of a cobordism $W$
between surfaces $\Sigma_{0},\Sigma_{1}$, of the kind considered
above, so $H_{1}(\Sigma_{0})\oplus H_{1}(\Sigma_{1})$ generates
$H_{1}(W)$. We let $U_{i}= H_{1}(\Sigma_{i})$ and let $\Gamma\subset
U_{0}\oplus U_{1}$ be the kernel of the inclusion map. Equivalently,
by Poincar\'e duality, we can take
$U_{i}=H^{1}(\Sigma_{i})$ and let $\Gamma$ be the image of the
restriction map, so under our hypotheses $\Gamma\cong H^{1}(W)$. We
use the integer lattices in all these spaces to fix the scalar
ambiguities, at least up to a sign, which we are not going to keep
track of. The conclusion is that we get a linear map
  $$\alpha_{W}\co \Lambda^{*}(H^{1}(\Sigma_{0}))
\rightarrow \Lambda^{*}(H^{1}(\Sigma_{1})). \tag 10$$ This can be
regarded as a prototype of the general construction recalled in
Section 2. Indeed we can obtain the map in a gauge theory framework
by considering moduli spaces of flat $S^{1}$ connections: the
exterior algebra then appears as the homology of the Jacobian torus
$H^{1}(\Sigma;\bold{R})/H^{1}(\Sigma;\bold{Z})$ parametrising flat
connections over $\Sigma$. Now consider how  the gradings work. In
general if $U_{i}$ has dimension $n_{i}$ and $\Gamma$ has dimension
$r$ the map $\vert \Gamma\vert$ increases degree by $n_{0}-r$. In
our case, Poincar\'e duality implies that
$$  \text{dim} (\Gamma)= \text{dim}(H^{1}(W)) =
 \frac{1}{2} (\text{dim} H^{1}(\Sigma_{0}) + \text{dim}
H^{1}(\Sigma_{1})), $$ and this means that $\alpha_{W}$ preserves the
{\it modified} grading of the exterior powers, so we have
$$ \alpha_{i,W}\co 
\Lambda_{(i)}(\Sigma_{0})\rightarrow\Lambda_{(i)}(\Sigma_{1}).
\tag11$$ As it stands, this notation is a little ambiguous because
$\Lambda_{(i)}(\Sigma)$ is used to denote the two (isomorphic)
exterior powers
$\Lambda^{g\pm i}$, and we have maps induced on each one. Let
$\omega_{i}$ be the intersection form on $H_{1}(\Sigma_{i})$, and
consider the symplectic form $(\omega_{0},-\omega_{1})$ on the
product. Then the subspace
$\Gamma_{W}$ is a {\it Lagrangian} subspace of the product with this
symplectic form---this is just the standard consequence of
Poincar\'e--Lefschetz duality for  $(W,\partial W)$. Then one can
readily show that the two maps associated to such a Lagrangian
subspace  are equal.

This construction behaves well with respect to the composition of
cobordisms. If
$W'$ is a cobordism from $\Sigma_{1}$ to $\Sigma_{2}$ and we form
the composite cobordism $W''$ by gluing $W$ to $W'$ along
$\Sigma_{1}$ and if $W''$ satisfies our homological condition
$H_{1}(W'',\partial W'')= \bold{Z}$ then the transversality
condition holds at $\Sigma_{1}$ and the composite subspace
$\Gamma''$ corresponds to $H^{1}(W'')$.

Now suppose that $\Sigma_{0}$ and $\Sigma_{1}$ are identified and we
glue the ends of $W$ to obtain a closed manifold $\overline{W}$ as
before. We also have a preferred element $\theta$ of
$H^{1}(\overline{W})$, the Poincar\'e dual of the image of the
boundary components, so we have an Alexander polynomial of the pair
$(\overline{W},\theta)$.
\proclaim{Proposition 12} The Alexander polyomial of
 $(\overline{W},\theta)$ is $a_{0}+ \sum_{j} a_{j}(t^{j}+ t^{-j})$
where
$$  a_{j} = (-1)^{j} \text{\rm Trace}(
\alpha_{j,W}\co\Lambda_{(j)}(\Sigma_{0})\rightarrow
\Lambda_{(j)}(\Sigma_{0})). $$
\endproclaim

Let $\tilde{W}$ be the infinite cyclic cover of $\overline{W}$
corresponding to
$\theta$. The Alexander polynomial is defined from the homology
 $H_{1}(\tilde{W})$,  regarded as a module over the ring
$\Lambda=\bold{Z}[t,t^{-1}]$ of Laurent series. Write $V_{\bold{Z}}$
for the homology group $H_{1}(\Sigma_{0})$ with integral
coefficients. We have a subgroup $\Gamma_{\bold{Z}}\subset
V_{\bold{Z}}\times V_{\bold{Z}}$, the integer lattice in $\Gamma$,
such that
$$   H_{1}(W;\bold{Z})= (V_{\bold{Z}}\oplus V_{\bold{Z}})/
\Gamma_{Z}. $$ Thinking of $\tilde{W}$ as constructed from a chain
of copies of $W$, glued along the boundaries, we see that, as an
abelian group
$$   H_{1}(\tilde{W})=\left(
 \dots V_{\bold{Z}}\oplus V_{\bold{Z}}\oplus
V_{\bold{Z}}\dots\right)/\sim, $$ where $\sim$ is the equivalence
relation generated by
$$     (\dots, 0,0,\sigma,0,0,\dots)\sim (\dots,0,0,0,\tau,0,\dots)
$$ if $(\sigma,\tau)\in \Gamma_{\bold{Z}}$. In other words
$$    H_{1}(\tilde{W}) = \Lambda \otimes_{\bold{Z}}
V_{\bold{Z}}/\sim, $$ where $\lambda\otimes \sigma\sim t \lambda
\otimes \tau $ for $(\sigma,\tau)\in \Gamma_{\bold{Z}}$. Explicitly,
this says that we can construct a square presentation matrix for
$H_{1}(\tilde{W})$ as  a $\Lambda$--module in the following way. Let
$(e_{j})$ be a standard basis for
$H_{1}(\Sigma_{0})$ and let  $(\gamma_{i})$ be a
 basis for $\Gamma_{\bold{Z}}$. We can write $\gamma_{i}=
(\sigma_{i},\tau_{i})
\in V_{\bold{Z}}\times V_{\bold{Z}}$, and express these in terms of
the basis:
$$
\sigma_{i}=\sum_{j} \sigma_{ij} e_{j} \ ,\
\ \tau_{i}= \sum_{j} \tau_{ij}e_{j} .$$ Then the $2g$ elements
$e_{i}$ can be thought of as  generators of $H_{1}(\tilde{W})$ over
$\Lambda$ with  $2g$ relations
  $$ \sum_{j} (\sigma_{ij}e_{j} - t \ \tau_{ij} e_{j}) = 0. $$
 So the presentation matrix is just $ (\sigma_{ij}- t \tau_{ij})$.
By definition the Alexander polynomial (up to a  unit in $\Lambda$)
is the determinant
$ \delta(t) = \det (\sigma_{ij}-t\tau_{ij})$. It is  now a
straightforward exercise in exterior algebra to complete the proof
of the proposition, up to  a multiplication by a unit $\pm t^{\mu}$,
 just working through the definitions in the Plucker construction.
Notice the particular  case when $\Gamma$ is the graph of a linear
map $f\co V\rightarrow V$. In this case we can choose bases so that
$(\sigma_{ij})$ is
 the identity matrix, and $\tau_{ij}$ is the matrix of $f$ in the
usual sense. The assertion that we want comes down to the familiar
fact that the coefficients of the characteristic polynomial are
$\pm \text{Tr} (\Lambda^{k} f)$.

With regard to the ambiguity in the definition of the Alexander
polynomial, recall  that this is reduced to $\pm 1$ by specifying
that $\Delta(t)=
\pm t^{\mu}\delta(t)$ should have $\Delta(t)=\Delta(t^{-1})$.  So to
complete the proof of Propostion 12,
 up to our overall $\pm1$ ambiguity, we need only use the fact noted
above that the traces of the $\alpha_{i,W}$ satisfy this symmetry
property.

We see now that the two formulae, for closed $3$--manifolds, can be
deduced if one proves the obvious  {\it relative} versions for
cobordisms. Suppose we have a theory in which we have fixed
isomorphisms, for all surfaces $\Sigma$,
$$   V(\Sigma) = \bigoplus_{j\geq 0}
 R_{j}\otimes \Lambda_{(j)}(\Sigma),\tag13 $$ where the $R_{i}$ are
{\it universal} vector spaces, independent of $\Sigma$.  Suppose we
can prove that for all cobordisms $W$ from $\Sigma_{0}$ to
$\Sigma_{1}$ the map $\rho_{W}\co  V(\Sigma_{0})\rightarrow
V(\Sigma_{1})$ is equal (under these fixed isomorphisms) to a direct
sum
$$       \rho_{W}= \bigoplus_{i}  1_{R_{j}}\otimes \alpha_{j, W}.
\tag14$$ Then in the case when $\Sigma_{0}=\Sigma_{1}$ the gluing
axioms and Proposition 12 imply that the numerical invariant of the
closed manifold $\overline{W}$ is
$$  \sum_{j} \text{dim}(R_{j})\  a_{j}. \tag15$$ In particular it is
then clear that
 the calculations of Section 3 will  imply the Casson and
Meng--Taubes formulae.  Now this relative version is both more
general and  easier to prove. The gluing axioms in our topological
field  theory and the functorial property of the Plucker
construction mean that if the statement is true for two composeable
cobordisms it is also true for the composite.
 We can decompose
 any cobordism $W$ into a composite of elementary cobordisms and it
suffices to prove the result for these. The homological condition
(3) means that we
 can choose these elementary cobordisms to have index $1$ or $2$;
that is, we reduce to considering the standard elementary cobordisms
$Z$ from $\Sigma^{g}$ to $\Sigma^{g+1}$ and
$Z'$ from $\Sigma^{g+1}$ to $\Sigma^{g}$, where  $\Sigma^{g}$ is a
standard surface of genus $g$. To sum up then we have  to establish,
 in the Seiberg--Witten and Casson theories, the existence of
isomorphisms (13) and the property (14) for these elementary
cobordisms.

\section{Naturality arguments}

We will now see that the proof can be completed without any further
geometric input, using  algebra and topology. In a nutshell, we
shall show that the invariants are determined by their formal
properties. We recall some representation theory of the (real)
symplectic group. The representations $\Lambda^{i}$ of
$Sp(2g,\bold{R})$ are not in general irreducible. Let
$L\co \Lambda^{i}\rightarrow \Lambda^{i+2}$  be the wedge product with
the symplectic form $\omega\in
\Lambda^{2}$. The {\it primitive subspace} $P^{i}$ is the kernel of
$$   L^{n-i+1}\co \Lambda^{i}\rightarrow \Lambda^{2n-i+2}. $$ The
decomposition of the exterior powers into irreducible representations
is, for $i\leq n$,
$$  \Lambda^{i} = P^{i} \oplus P^{i-2} \oplus \dots , \tag16$$ [6],
where we use the maps $L^{j}$ to map $P^{i-2j}$ into $\Lambda^{i}$.
Thus
$$   P^{i} = \Lambda^{i}-\Lambda^{i-2}, $$ as virtual
representations. We shall also need to
 use the fact that the $P^{i}$ are irreducible representations of
the discrete subgroup $Sp(2n,\bold{Z})\subset Sp(2n,\bold{R})$, and
of the finite-index subgroups of $Sp(2n,\bold{Z})$ associated to
spin structures.    This refinement
 follows because these subgroups are Zariski-dense in
$Sp(2n,\bold{R})$.

 We now modify the grading of the primitive spaces, writing
$P_{(i)}= P^{n-i}$, so
$$    \Lambda_{(i)} = P_{(i)} + P_{(i+2)} + \dots, . \tag17$$ We
associate to a surface $\Sigma$ vector spaces
$$P_{(i)}(\Sigma)= P_{(i)}(H_{1}(\Sigma)),$$ which give irreducible
representations of the mapping class group. Now suppose we have a
theory in which we know that
$$ V(\Sigma) = \sum \mu_{j} \Lambda_{(j)}, $$ as representations of
the mapping class group, for certain multiplicities
$\mu_{i}$.  Suppose the theory is $\bold{Z/2}$--graded, with the
even/odd part made of the $\Lambda_{(\text{even})}$,\break
$\Lambda_{(\text{odd})}$ respectively. Suppose the multiplicities
$\mu_{j}$ are {\it universal}, independent of the genus. All of
these conditions hold in our two geometric examples. Then we can
write
$$  V(\Sigma) = \sum_{j} \nu_{j} P_{(j)}, $$ where
$$\nu_{0}=\mu_{0}, \nu_{1}=\mu_{1}, \nu_{2}=\mu_{0}+\mu_{2}, \dots$$
Define vector spaces $Q_{j}$ for $0\leq j\leq g(\Sigma)$ by
$$   Q_{j}(\Sigma) = \text{Hom}_{MC(\Sigma)} (P_{(j)}(\Sigma),
V(\Sigma))
$$ where $\text{Hom}_{MC(\Sigma)}$ denotes the equivariant
homomorphisms under the action of the mapping class group. So the
$Q_{j}$
 are vector spaces canonically associated to a surface and by
construction the mapping class group acts trivially on them.
 On the other hand we have {\it canonical} isomorphisms
$$   V(\Sigma) = \bigoplus_{j} Q_{j}(\Sigma)\otimes P_{(j)}(\Sigma),
\tag18$$ since the $P_{(j)}$ are distinct irreducible
representations  of $Sp(2n,\bold{R})$ and hence of the mapping class
group. (Here we use the fact that the mapping class group maps onto
$Sp(2n\bold{Z})$.)  Notice that the dimensions of the
$Q_{j}(\Sigma)$ are the integers
$\nu_{j}$ determined by the $\mu_{j}$.  and hence are universal,
independent of the genus.

At this stage we bring in the following lemma, whose proof we leave
as an exercise for the reader. (For our main  proof we only need a
special case of  this---for the maps induced by elementary
cobordisms---where the exercise is more straightforward.)

\proclaim{Lemma 19}
   Suppose $(U_{0},\omega_{0}), (U_{1},\omega_{1})$
 are symplectic vector spaces and $\Gamma\subset U_{0}\times U_{1}
$ is a Lagrangian subspace with respect to the    symplectic form
$(\omega_{0}, -\omega_{1})$ on $U_{0}\times U_{1}$.  Then the
Plucker map $\vert \Gamma \vert$ maps $P_{(j)}(U_{0})$ to
$P_{(j)}(U_{1})$.
\endproclaim

It follows that a cobordism $W$ defines maps
$$  \beta_{j,W}\co  P_{(j)}(\Sigma_{0})\rightarrow P_{(j)}(\Sigma_{1}),
\tag20$$ since the corresponding subspace $\Gamma_{W}\subset
H_{1}(\Sigma_{0})\times  H_{1}(\Sigma_{1})= H_{1}(\partial W)$ is
Lagrangian. Consider the standard cobordism $Z$ from $\Sigma^{g}$ to
$\Sigma^{g+1}$. This defines
$$  \rho_{Z}\co  \bigoplus_{j} Q_{j}(\Sigma^{g})\otimes
P_{(j)}(\Sigma^{g})
\rightarrow \bigoplus_{j} Q_{j}(\Sigma^{g+1}) \otimes
P_{(j)}(\Sigma^{g+1}),
 \tag21$$ where we have used the identifications (18) on the two
ends of the cobordism.

\proclaim{Lemma 22} There are  linear maps $\rho_{j}\co 
Q_{j}(\Sigma^{g})\rightarrow  Q_{j}(\Sigma^{g+1})$ such that
$$\rho_{Z}= \bigoplus_{j} \rho_{j}\otimes \beta_{j,Z}. $$
\rm

We may construct the elementary cobordism $Z$ in a canonical way
starting with  the surface $\Sigma^{g+1}$ and a non-separating
embedded circle $\delta
 \subset \Sigma^{g+1}$.
 The homology group $H_{1}(\Sigma^{g})$ can be identified with the
quotient  $[\delta]^{\perp}/[\delta]$ of the annihilator
$[\delta]^{\perp}\subset H_{1}(\Sigma^{g+1})$ of $[\delta]\in
H_{1}(\Sigma^{g+1})$ with respect to the intersection form. Any
diffeomorphism of $\Sigma^{g+1}$ that fixes $\delta$ induces a
diffeomorphism of $Z$. Let $G$ be the group of  automorphisms of
$H_{1}(\Sigma^{g+1};\bold{Z})$ that preserve the intersection form
and fix $[\delta]$. We get a natural action of $G$ on
$H_{1}(\Sigma^{g})= [\delta]^{\perp}/\delta$. Any element of $G$ can
be lifted to a diffeomorphism of $\Sigma^{g+1}$, fixing $\delta$,
and hence to a diffeomorphism of $Z$. It follows from the naturality
properties of our theory that $\rho_{Z}$ intertwines the
$G$--actions on the spaces $\bigoplus Q_{j}\otimes P_{(j)}$ at the
two ends, and these actions are induced from the actions on the
$P_{(j)}$ since the mapping class groups act trivially on the
$Q_{j}$. Now Lemma 22 follows from the $\bold{Z}/2$--grading in the
theory
 and the next, purely algebraic, Lemma.
\proclaim{Lemma 23} If $j+j'$ is even the only non-zero
$G$--equivariant maps from $P_{(j)}(\Sigma^{g})$ to $
P_{(j')}(\Sigma^{g+1})$ are multiples of $\beta_{j,Z}$, in the case
when $j=j'$.

\endproclaim

To see this we first choose another class $[\epsilon]\in
H_{1}(\Sigma^{g+1})$ with $[\epsilon].[\delta]= 1$. This choice
defines an isomorphism
$$   H_{1}(\Sigma^{g+1}) = H_{1}(\Sigma^{g}) \oplus \langle
[\epsilon],[\delta]\rangle , $$ and hence a copy of
$Sp(2g,\bold{Z})$ in $ G$. It is easy to verify that we have then
have natural (ie $Sp(2g,\bold{Z})$--invariant)  isomorphisms
$$   P_{(j)}(\Sigma^{g+1}) =
 P_{(j+1)}(\Sigma^{g}) \oplus P_{(j)}(\Sigma^{g})\otimes \langle
[\epsilon, \delta]
\oplus P_{(j-1)}(\Sigma^{g}). $$ Using the fact that the $P_{(j)}$
are distinct irreducible representations we see that the only
$Sp(2g,\bold{Z})$--equivariant maps are the maps in this
decomposition   from $P_{(j)}(\Sigma^{g})$ to
$P_{(j)}(\Sigma^{g})\otimes
\langle [\epsilon],[\delta]\rangle$ of the form $1\otimes \lambda$,
for
$\lambda\in \langle [\delta],[\epsilon]\rangle$. But the only
classes in
$\langle [\delta],[\epsilon]\rangle$ which are fixed by the larger
group $G$ are the multiples of $\delta$, so the only $G$--equivariant
maps are multiples of $1\otimes \delta$ which are just the multiples
of $\beta_{j,Z}$.

Note that it follows from the definitions, and the trivial action on
the
$Q_{j}$,  that the maps $\rho_{j}$ are universal in the following
sense. If we start with surfaces $\Sigma^{g}, \Sigma^{g+1}$ and
choose any standard cobordism between them we get the same maps
$\rho_{j}$.

\proclaim{Lemma 24} The maps $\rho_{j}$ are isomorphisms, for $0\leq
j\leq g$. \endproclaim

We may carry through the entire discussion above for the other
elementary cobordism
$Z'$ from $\Sigma^{g+1}$ to $\Sigma^{g}$. We get
$$ \rho_{Z'} = \bigoplus \rho'_{j} \otimes \beta_{j,Z'}, $$ for
universal maps $\rho'_{j}$. Consider the trivial cobordism
$\Sigma^{g}\times[0,1]$ and write this as a composite of elementary
cobordisms
$Z, Z'$. That is, we introduce a pair of cancelling handles in the
product.
 The map in our theory induced by the trivial cobordism is the
identity map, and it follows that $\rho'_{j}\circ \rho_{j}$ is the
identity for each $j$. So  the $\rho_{j}$ are injective. But we know,
 by our original hypothesis, that for $0\leq j\leq g$
$$\text{dim}\ Q_{j}(\Sigma^{g})= \text{dim}\ Q_{j}(\Sigma^{g+1}), $$
so the $\rho_{j}$ are isomorphisms, with inverses $\rho_{j}'$.

Now  we use the universal maps $\rho_{j}$ to {\it identify} the
spaces
$Q_{j}(\Sigma)$ for all surfaces $\Sigma$ (of genus $g\geq j$),
 and hence define universal vector spaces
$Q_{j}$.   The point is that this is unambigous since the mapping
class groups acts trivially on the
$Q_{j}(\Sigma)$ so it does not matter which cobordisms we use to
induce the identifications.  Then we have fixed isomorphisms
$$ V(\Sigma) = \bigoplus_{j} Q_{j}\otimes P_{(j)}(\Sigma), $$ which,
by construction,  commute with the maps induced by elementary
cobordisms, and hence by all cobordisms. Finally we can get back to
the exterior powers,  although this is not really necessary, since
we can go straight from the primitive spaces to the Alexander
polynomial.  We know that the dimension of $Q_{j}$ is
$\mu_{j}+\mu_{j-2}+ \dots $, so we may choose {\it arbitrary}
embeddings
$$ \align  Q_{0}&\hookrightarrow Q_{2}\hookrightarrow Q_{4} \dots \\
           Q_{1}&\hookrightarrow Q_{3}\hookrightarrow Q_{5}\dots
\endalign$$ and a sequence of complementary subspaces
$$     Q_{i}= Q_{i-2}\oplus R_{i}, $$ so that $$   Q_{i}=
R_{i}\oplus R_{i-2}\oplus R_{i-4} \dots, $$ and $\text{dim}
R_{i}=\mu_{i}$. Then we have fixed isomorphisms
$V(\Sigma)=\bigoplus R_{j}\otimes \Lambda_{(j)}(\Sigma)$, which
commute with the maps in the theory, as required.

\section{Connections with work of Thaddeus}

The main point of this article has been the fact that the homology of
the flat connection moduli spaces, and of the symmetric products,
can be described in terms of the exterior powers. It is natural to
look  for relations between
 these different geometric problems.  In his paper [12], Thaddeus
used an algebro--geometric construction to relate these
objects---regarding the moduli space $M(\Sigma)$ as a moduli space
of rank
$2$ holomorphic vector bundles. His idea is to study a moduli space
of
 \lq\lq stable pairs''
$(E,\phi)$, where $E$ is a bundle and $\phi$ is a holomorphic
section of
$E$. The notion of stability depends on a real parameter
 $\tau\in (0,\infty)$, so there is a family of moduli spaces
$\Cal{M}_{\tau}$. When $\tau$ is small this fibres over the moduli
space of vector bundles with fibre a complex projective space, and
when $\tau$ is large
$\Cal{M}_{\tau}$ is simply a projective space. In between there are
a finite number of exceptional values of $\tau$ where
$\Cal{M}_{\tau}$ changes by a \lq\lq flip'' (or \lq\lq complex
surgery'') which can be described in terms of the symmetric
products.  So far in this paper we have avoided the case $d=0$ in
the Seiberg--Witten theory. In this case there are reducible
solutions of the equations over a surface,
 which causes some complications. However one can make a
perturbation of the equation which removes the reducible solutions,
and the relevant moduli space is then the  symmetric product
$\sym^{g-1}(\Sigma)$. Thus we can  extend our theory with vector
spaces
$$   V_{0}= \Lambda_{(1)}\oplus 2 \Lambda_{(2)} \oplus \dots . $$
Similarly, if we use perturbations, we obtain a theory for {\it
negative} values of $d$ which brings in the higher symmetric
products.
 (This is not the same as the unperturbed theory, in which there is
a symmetry between $d$ and $-d$.) Thus
$$   V_{-d} = \Lambda_{-(d-1)} \oplus 2\Lambda_{-(d-2)}\oplus \dots,
$$ but using the symmetry $\Lambda_{(k)}=\Lambda_{(-k)}$ we get:
$$    V_{-d} = V_{d} \oplus d T, \tag 25$$ where
$$  T= \Lambda_{0} \oplus 2\Lambda_{1} \oplus \Lambda_{2} \dots\ . $$
(This corresponds to the homology of the Jacobian torus.) Now
Thaddeus' relation between the homology of the different spaces can
be written
$$ (2g+1)  H_{*}(M(\Sigma)) = \sum_{j=0}^{2g-1}
 (5g-2-3j) H_{*}(\sym^{j}(\Sigma))= \sum_{k=-d}^{k=d-1}
(2g+1+3d)V_{d}. $$ Using (25), the last sum yields the relation
$$   V^{C} = V_{0} + 2V_{1} + 2 V_{2} + 2V_{3} + \dots , $$ which we
can of course obtain equally well directly from (5), (9).
 This can be seen as a
 universal formula, relating the
$(2+1)$--dimensional Casson theory and the $(2+1)$--dimensional
Seiberg--Witten
 theory. These are both reductions of $(3+1)$ dimensional
theories---the instanton theory in the Casson case---so the ideas
should have some bearing on the problem of understanding the
relation between these latter theories. Of course there are obvious
similarities between Thaddeus' technique and those used in  the
programme of Pidstragatch and Tyurin in
$4$--dimensions. It would be interesting to obtain the relation
between the theories directly using the Pidstragatch and Tyurin
method.

\np

\Addresses\recd

\enddocument